\documentclass[french,english,reqno,12pt]{amsart}
\usepackage{amssymb,amsmath}
\usepackage{euscript}
\usepackage{ifthen}
\usepackage{array}
\usepackage[all]{xy}
\usepackage{babel}
\theoremstyle{plain}
\newtheorem{theorem}{Theorem}[section]
\newtheorem{proposition}[theorem]{Proposition}
\newtheorem{lemma}[theorem]{Lemma}
\newtheorem{corollary}[theorem]{Corollary}

\theoremstyle{definition}

\newcommand{\M}{\EuScript{M}}

\newcommand{\QC}{\widehat{\mathcal{Q}}}
\renewcommand{\S}{S}
\newcommand{\calS}{\mathcal{S}}
\newcommand{\calQ}{\mathcal{Q}}
\newcommand{\bfX}{{\bf X}}
\newcommand{\bfY}{{\bf Y}}
\newcommand{\bfe}{{\bf e}}
\newcommand{\bfone}{\mathbf{1}}

\newcommand{\field}{\Bbbk}
\newcommand{\Rep}{\operatorname{Rep}}
\newcommand{\Res}{\operatorname{Res}}
\newcommand{\Ind}{\operatorname{Ind}}
\newcommand{\End}{\operatorname{End}}
\newcommand{\EndGr}{\operatorname{end}}
\newcommand{\Sh}{\operatorname{Sh}}
\newcommand{\GL}{\operatorname{GL}}
\newcommand{\Des}{\operatorname{Des}}

\DeclareMathOperator{\smashprod}{\#}

\newcounter{xx}
\newcommand{\ver}{\ifthenelse{\value{xx} = 0}{\vrule}{}}

\newcommand{\xyinc}{\ar@{^{(}->}}

\newcommand{\onto}{\twoheadrightarrow} 

\newcommand{\inc}{\hookrightarrow}

\newcommand{\chir}{\raisebox{2pt}{$\chi$}}
\newcommand{\comp}{\vDash}              %for compositions (of n)
\newcommand{\parti}{\vdash}              %for partitions (of n)

\begin{document}
\selectlanguage{english}

\title[The smash product of symmetric functions]{The smash product of symmetric functions.\\ Extended abstract}
 \author[M.~Aguiar, W. Ferrer, and W. Moreira]{Marcelo Aguiar, Walter Ferrer, and Walter Moreira}

 \address{Department of Mathematics\\ Texas A\&M University\\
 College Station, TX 77843, USA}
 \email{maguiar@math.tamu.edu}
\urladdr{http://www.math.tamu.edu/$\sim$maguiar}

 \address{Centro de Matem\'atica\\  Universidad de la Rep\'ublica\\
Igu\'a 4225 esq. Mataojo\\
 Montevideo 11400, Uruguay}
 \email{wrferrer@cmat.edu.uy}
\urladdr{http://www.cmat.edu.uy/$\sim$wrferrer}

\address{Department of Mathematics\\ Texas A\&M University\\
 College Station, TX 77843, USA}
 \email{wmoreira@math.tamu.edu}
\urladdr{http://www.math.tamu.edu/$\sim$wmoreira}

\thanks{Aguiar supported in part by NSF grant DMS-0302423}
\thanks{The authors thank Arun Ram and Frank Sottile for interesting conversations and suggestions}

\date{November 21, 2004}

%\begin{abstract}
%\end{abstract}

\maketitle

\begin{abstract}
  We construct a new operation among representations of the
  symmetric group that interpolates between the classical {\em internal} and
 {\em external} products, which are defined in terms of tensor product
 and induction of representations.  Following Malvenuto and
 Reutenauer, we pass from symmetric functions to non-commutative
 symmetric functions and from there to the algebra of permutations in
 order to relate the internal and external products to
  the {\em composition} and {\em convolution} of linear endomorphisms of the tensor algebra.
  The new product we construct corresponds to the {\em smash} product 
   of endomorphisms of the tensor algebra.  For symmetric functions, the
   smash product is given by a construction which combines induction and
   restriction of representations. 
   For non-commutative symmetric functions, the structure constants of
   the smash product are given by an explicit combinatorial rule which
   extends a well-known
 result of Garsia, Remmel,  Reutenauer, and  Solomon for the descent algebra.
  We describe the dual operation among quasi-symmetric functions in terms of alphabets.
\end{abstract}

\begin{otherlanguage}{francais}
\begin{abstract}
 Nous construisons une nouvelle op\'eration parmi les
 repr\'esentations du groupe sym\'etrique qui interpole entre les
 produits {\em interne} et {\em externe}. Ces derniers sont d\'efinis
 en termes du produit tensoriel et de l'induction des
 repr\'esentations. D'apr\`es Malvenuto et Reutenauer, nous passons
 des fonctions sym\'etriques aux fonctions sym\'etriques non
 commutatives et \`a l'alg\`ebre des permutations afin de rapporter
 les produits internes et externes \`a la composition et \`a la
 convolution d'endomorphismes lin\'eaires de l'alg\`ebre
 tensorielle. Le nouveau produit correspond au produit {\em smash}
 d'endomorphismes de l'alg\`ebre tensorielle. Pour les fonctions
 sym\'etriques, le produit smash est donn\'e par une construction qui
 combine l'induction et la restriction de repr\'esentations. Pour les
 fonctions sym\'etriques non commutatives, les constantes de structure
 du produit smash sont donn\'ees par une r\`egle combinatoire
 explicite qui prolonge un r\'esultat bien connu de Garsia, Remmel,
 Reutenauer et Solomon pour l'alg\`ebre de descentes. Nous d\'ecrivons
 l'op\'eration duale au niveau des fonctions quasi-sym\'etriques en
 termes d'alphabets.  \end{abstract} \end{otherlanguage}

\selectlanguage{english}

\reversemarginpar

\section*{Introduction}

Our goal is to introduce a new operation among symmetric functions which interpolates between the classical internal and external products.
This operation is best understood by considering not only symmetric functions but also three other algebras, related by means of
the following fundamental commutative diagram:
\begin{equation}
\begin{gathered}\label{E:fund-diagram}
\xymatrix@R=15pt@C=15pt{
{\;\Sigma\;}\xyinc[r]\ar@{>>}[d]  &{\;\calS\;}\ar@{>>}[d] \\ 
{\;\Lambda\;}\xyinc[r]    &{\;\calQ\;} }
\end{gathered}
\end{equation}
$\Lambda$ is the algebra of symmetric functions~\cite{Mac,Sag,Sta}, $\calQ$ is the
algebra of quasi-symmetric functions~\cite{Ges}, $\Sigma$ is the algebra of non-commutative symmetric functions~\cite{GKLLRT}, and $\calS$ is the algebra of permutations~\cite{MR}. All four are graded connected Hopf algebras.
The definitions of these Hopf algebras, as well as the maps that relate them,
are reviewed in this paper.

We work over a field $\field$ of characteristic $0$.

Let $V$ be a representation of the symmetric group $\S_n$. The {\em Frobenius characteristic} of $V$ is the symmetric function
\[\sum_{\lambda\parti n}
\frac{\chir_V(\lambda)}{z_\lambda}\sum_{i_1,\ldots,i_r}x_{i_1}^{\ell_1}\cdots
x_{i_r}^{\ell_r}\,.\]
The sum is over all partitions  $\lambda=(\ell_1\geq\cdots\geq\ell_r)$ of $n$,  $z_\lambda$ is the order of the stabilizer of the conjugacy class of permutations of cycle-type $\lambda$, and
$\chir_V(\lambda)$ is the character of $V$ evaluated on any such permutation.
 This association allows us to identify the algebra of symmetric functions with the
 direct sum of the Grothendieck groups of $\S_n$:
 \[\Lambda =\bigoplus_{n\ge 0} \Rep(\S_n)\,.\]
Let $V$ be a representation of $\S_p$ and $W$ a representation of $\S_q$.
The {\em external} product of $V$ and $W$ is the representation
\[V\ast W:=\Ind_{\S_p\times \S_q}^{S_{p+q}}(V\otimes W)\]
of $\S_{p+q}$. This operation corresponds to the product of power series under the Frobenius characteristic. If $V$ and $W$ are representations of $\S_n$, their {\em internal} product is the diagonal representation of $\S_n$ on their tensor product:
\[V\circ W:= \Res^{\S_n\times \S_n}_{\S_n}(V\otimes W)\,.\]
Explicit expressions for these products on the basis of irreducible representations
(Schur functions) are of central interest in the theory of symmetric functions.
While a complete solution for the case of the external product is known
(the {\em Littlewood-Richardson rule}), only partial answers are known for the case of the internal product (the {\em Kronecker problem}).

In Section~\ref{S:rep} we introduce a new product $V\smashprod W$ between representations of $\S_p$ and $\S_q$ which contains the internal and external products as the terms of extreme degrees, as well as additional terms, and which is still associative. We say that this {\em smash product} of representations interpolates between the internal and external products. 
For example, the smash product of the complete symmetric
functions $h_{(2,1)}$ and $h_3$ is
\begin{equation*}
h_{(2,1)}\smashprod h_3 = h_{(2,1)}+ h_{(1,1,1,1)}+ h_{(2,1,1)} +
h_{(2,2,1)}+h_{(2,1,1,1)} + h_{(3,2,1)},
\end{equation*}
where the external product is recognized in the last term and the
internal product in the first one, together with additional terms of
degrees four and five.

The existence of this operation  poses the problem of finding an explicit description for its structure constants on the basis of Schur functions.
The answer would contain as extreme cases the Littlewood-Richardson rule and
(a still unknown) rule for the Kronecker coefficients.

The smash product arises from a  construction in the theory of Hopf algebras which is a simple generalization  of the notion of semidirect product of groups.
Let $H$ be a Hopf algebra. Considering the action by
translations of $H$ on its dual $H^*$  leads to an associative operation on the space $\End(H)$ of linear endomorphisms of $H$ (Section~\ref{S:end-per}).
 Let $H=\bigoplus_{n\geq 0} H_n$ be a graded connected Hopf algebra. In this general setting, the smash product of two endomorphisms $f:H_p\to H_p$ and $g:H_q\to H_q$ is a sum of various endomorphisms $H_n\to H_n$, with $\max(p,q)\leq n\leq p+q$. The endomorphism corresponding to $n=p+q$ is the familiar {\em convolution} of $f$ and $g$, while that one corresponding to $n=p=q$ is simply the composition of $g$ and $f$.

The connection to symmetric functions is made through diagram~\eqref{E:fund-diagram}, starting from the opposite vertex.  We take $H=T(V)$, the tensor algebra of a vector space.
First of all, Schur-Weyl duality allows us to restrict the smash product of endomorphisms of $T(V)$ to the direct sum of the symmetric group algebras
 \[\calS=\bigoplus_{n\ge 0}\field\S_n\]
(Corollary~\ref{C:smash-per}). Here  $\S_n$ is viewed as the endomorphisms of $V^{\otimes n}$ which permute the coordinates. The convolution product on $\calS$ is the product of Malvenuto and Reutenauer,
 while the composition product is simply the group algebra product.
Let
 \[\Sigma=\bigoplus_{n\ge 0}\Sigma_n\]
 be the direct sum of {\em Solomon's descent algebras}. A result of Garsia and Reutenauer which characterizes the elements of $\Sigma$ in terms of the action on the tensor algebra allows us to show that the smash product also restricts to $\Sigma$~(Theorem~\ref{T:smash-sigma}). In particular, we recover the classical result of Solomon that each $\Sigma_n$ is
 a subalgebra of the symmetric group algebra $\field\S_n$, and the fact that
 $\Sigma$ is closed under the convolution of permutations. Endowed with the latter, $\Sigma$ is the algebra of non-commutative symmetric functions.
 
 In  Theorem~\ref{T:smash-X} we  provide a combinatorial rule which describes the structure constants of the smash product of $\Sigma$ on the basis of shuffles (the complete non-commutative symmetric functions). This rule interpolates nicely between the well-known rule of Garsia, Remmel, Reutenauer, and Solomon for the descent algebras and the concatenation rule for the product of complete non-commutative symmetric functions.
 
We arrive at the smash product of symmetric functions by showing that the smash product of $\Sigma$ corresponds to the smash product of $\Lambda$ via Solomon's epimorphism $\Sigma\onto\Lambda$ (Theorem~\ref{T:sigma-lambda}).

In Section~\ref{S:quasi} we provide a description for the coproduct of quasi-symmetric functions corresponding to the smash product of non-commutative symmetric functions by duality. The description is in terms of alphabets (Theorem~\ref{T:coprod-Q}).

We also discuss the Hopf algebra structure that accompanies  the smash product on $\Sigma$, $\Lambda$ and $\calQ$.

\section{The smash product of representations of the symmetric group}
\label{S:rep}

Let $p,q,n$ be non-negative integers with $\max(p,q)\leq n\leq p+q$. 
Given permutations $\sigma\in \S_p$ and $\tau\in \S_q$, let $\sigma\times\tau\in \S_{p+q}$ be
\begin{equation}\label{E:times}
(\sigma\times\tau)(i)=\begin{cases} \sigma(i) & \text{ if $1\leq i\leq p$,}\\
\tau(i-p)+p & \text{ if $p+1\leq i\leq p+q$.} \end{cases}
\end{equation}

Let $\S_p\times_n \S_q:=\S_{n-q} \times\S_{p+q-n} \times\S_{n-p}$.
 Consider the embeddings
%\begin{equation*}
%\xymatrix@R=15pt@!C=8pt{ \S_n & & \S_p \times \S_q  \\
%& \S_p \times_n \S_q \ar@{_{(}->}[ul] \ar@{^{(}->}[ur] }
%\end{equation*}
\begin{align}
 \S_p \times_n \S_q \inc \S_n & \quad (\sigma,\rho,\tau)\mapsto \sigma\times\rho\times\tau\,,  \label{E:embed-n}\\
 \S_p \times_n \S_q \inc \S_p\times\S_q & \quad (\sigma,\rho,\tau)\mapsto (\sigma\times\rho,\rho\times\tau)\,.  \label{E:embed-pq}
\end{align}

The \textit{smash product} of a representation $V$ of $\S_p$ and a representation $W$ of $\S_q$ is
\begin{equation}\label{E:smash-rep}
V \smashprod W := \bigoplus_{n=\max(p,q)}^{p+q} \Ind_{\S_p\times_n
  \S_q}^{\S_n}\Res_{\S_p\times_n\S_q}^{\S_p\times\S_q} (V\otimes W)\,.
\end{equation}
This is an element of the direct sum of the Grothendieck groups of the symmetric groups
\[\Lambda =\bigoplus_{n\ge 0} \Rep(\S_n)\,.\]
We work over a field $\field$ of characteristic $0$; $\Lambda$ is a graded vector space over $\field$.

Let $(V\smashprod W)_n$ denote the component of degree $n$ in~\eqref{E:smash-rep}.
Consider the top component ($n=p+q$). In this case 
embedding~\eqref{E:embed-pq} is the identity and~\eqref{E:embed-n}
is the standard parabolic embedding $\S_p \times \S_q \inc \S_{p+q}$. We thus get
\begin{equation*}
(V\smashprod W)_{p+q}= \Ind_{\S_p\times \S_q}^{\S_{p+q}}(V\otimes
W),
\end{equation*}
which is the usual external product of representations~\cite{Gei,Zel}.

On the  other hand, when $n=p=q$, embedding~\eqref{E:embed-n} is the identity and~\eqref{E:embed-n} is the diagonal embedding $\S_n \inc \S_n\times\S_n$.
Therefore,
\begin{equation*}
(V\smashprod W)_n = \Res^{\S_n\times \S_n}_{\S_n}(V\otimes W),
\end{equation*}
 the internal product of representations (also known as Kronecker's
product)~\cite{Gei,Zel}. 

The smash product contains terms of intermediate degrees between
$\max(p,q)$ to $p+q$; in this sense it
``interpolates'' between the internal and external products.  It is a remarkable
fact that, as the internal and external products, the smash product is associative,
and can be lifted  to other settings (non-commutative symmetric functions,
permutations, and dually, quasi-symmetric functions).

Let $\alpha=(a_1,\ldots,a_r)$ be a composition of $n$ and
\[\S_\alpha:=\S_{a_1}\times\cdots\times\S_{a_r}\,.\]
We view $\S_\alpha$ as a subgroup of $\S_n$ by iterating~\eqref{E:times}. 
Let $h_\alpha$ denote the permutation representation of $\S_n$ corresponding to the action by multiplication on the quotient $\S_n/\S_\alpha$. The isomorphism class of $h_\alpha$ does not depend on the order of the parts of $\alpha$. As $\alpha$ runs over the set of partitions of $n$, the representations $h_\alpha$ form a linear basis of $\Rep(\S_n)$.

We provide an explicit description for the smash product on this
basis. Let  $\alpha=(a_1,\ldots,a_r)\comp p$  and  $\beta=(b_1,\ldots,b_r)\comp q$ be two compositions and $n$ an integer  with
$\max(p,q)\leq n\leq p+q$. Let $a_0=n-p$, $b_0=n-q$, and
 let $\M_{\alpha,\beta}^n$ be the set
of all $(s+1)\times(r+1)$-matrices 
\[M=(m_{ij})_{0\leq i\leq s,\,0\leq j\leq r}\]
 with non-negative integer entries and such that 
 \begin{itemize}
\item the sequence of column sums  is $(a_0,a_1,\ldots,a_r)$,
\item the sequence of row sums  is $(b_0,b_1,\ldots,b_s)$,
\item the first entry is $m_{00}=0$.
 \end{itemize}
We illustrate these conditions as follows.
\setcounter{xx}{0}
\begin{equation*}
    \begin{array}{cccc!{\ver}c}
      0 & m_{01} & \cdots & m_{0r} & n-q \\ 
      m_{10} & m_{11} & \cdots & m_{1r} & b_1 \\
      \vdots & \vdots & \ddots & \vdots & \vdots \\
      m_{s0} & m_{s1} & \cdots & m_{sr} & b_s \\ 
      \cline{1-4} 
      n-p & a_1  & \cdots &\setcounter{xx}{1} a_r & \\
    \end{array}
  \end{equation*}
 Let $p(M)$ be the partition of $n$ whose parts are the non-zero $m_{ij}$. 

\begin{theorem} \label{T:smash-rep}
 \begin{equation*}
 h_\alpha \smashprod h_\beta = \bigoplus_{n=\max(p,q)}^{p+q}
 \bigoplus_{M \in \M_{\alpha,\beta}^n} h_{p(M)}\,.
\end{equation*}
\end{theorem}

The proof of Theorem~\ref{T:smash-rep} follows from the Mackey rule from representation theory.

The space $\Lambda$ can be equipped with a  coproduct~\cite{Gei,Zel}.
 View  
 \[\Rep(\S_p \times \S_q) =  \Rep(\S_p) \otimes \Rep(\S_q)\,.\]
Given $V \in \Rep(\S_n)$ define  
\begin{equation}\label{E:coprod-rep}
 \Delta(V)= \sum_{p+q=n}  \Res^{\S_n}_{\S_p \times \S_q}(V)\,,
\end{equation}
where the restriction is along the parabolic embedding~\eqref{E:times}.
\begin{theorem}\label{T:Hopf-rep}
The space $\Lambda$ endowed with the smash product~\eqref{E:smash-rep} and the coproduct~\eqref{E:coprod-rep} is a  connected Hopf algebra. It is
commutative and cocommutative.
\end{theorem}

This is deduced from a similar statement
for the space of non-commutative symmetric functions $\Sigma$ (Theorem~\ref{T:Hopf-sigma}) via the
canonical surjection $\Sigma\onto\Lambda$.
This is done in Section~\ref{S:sym}.

%The compatibility of the smash product and the coproduct is proved
%using the explicit formula for the coefficients, which allow us to
%easily establish a bijection between the terms of of
%$\Delta(X_\alpha\smashprod X_\beta)$ and
%$\Delta(X_\alpha)\smashprod\Delta(X_\beta)$.
%\begin{comm} Should we mention here also the smash coproduct?
%\end{comm}

\section{The smash product of endomorphisms and permutations} \label{S:end-per}

Let $H$ be an arbitrary Hopf algebra, let $m:H\otimes H\to H$ be the product and $\Delta:H\otimes H\otimes H$ be the coproduct. The space $\End(H)$ of linear
endomorphisms of $H$ carries several associative products. Let $f,g\in\End(H)$.
Composition and convolution are respectively defined by the diagrams
\begin{equation}\label{E:com-con}
\begin{gathered}
\xymatrix@R=15pt@C=8pt{
&  {H}\ar[rd]^{f} &\\
{H}\ar[ru]^{g}\ar[rr]_{f\circ g} & & {H}  } \qquad
\xymatrix@R=15pt@C=18pt{
{H\otimes H}\ar[r]^{f\otimes g} & {H\otimes H}\ar[d]^{m}\\
{H}\ar[u]^{\Delta}\ar[r]_{f\ast g} & {H}}
\end{gathered}
\end{equation}
The smash product of endomorphisms is  defined by the diagram
 \begin{equation}\label{E:smash-end}
\begin{gathered}
\xymatrix@R=15pt@C=8pt{
&  {H^{\otimes 3}\ }\ar[rr]^{\mathrm{cyclic}} &&{\ H^{\otimes 3}}\ar[rd]^{1\otimes m}\\
{H^{\otimes 2}}\ar[ru]^{\Delta\otimes 1} & && & {H^{\otimes 2}}\ar[d]^{1\otimes g}\\
H^{\otimes 2}\ar[u]^{f\otimes 1}& && & {H^{\otimes 2}}\ar[ld]^{m}\\
& H\ar[lu]^{\Delta}\ar[rr]_{{f\smashprod g}} && H }
\end{gathered}
\end{equation}
where the map $\mathrm{cyclic}:H^{\otimes 3}\to H^{\otimes 3}$ is
$x\otimes y\otimes z\mapsto y\otimes z\otimes x$. 
Associativity follows from the Hopf algebra axioms.

The smash product is often defined in a different setting~\cite{Mon}:
given a Hopf algebra $H$ and an $H$-module-algebra $A$, the smash product
is an operation on the space $A\otimes H$ defined by
\begin{equation}\label{E:smash-usual}
(a\otimes h)\smash(b\otimes k):=\sum a(h_1\cdot b)\otimes h_2k\,.
\end{equation}
If $A=H^*$ and $H$ acts on $A$ by translation then~\eqref{E:smash-usual} corresponds to~\eqref{E:smash-end} via the canonical inclusion $H^*\otimes H\inc\End(H)$. Note that we have not made any finite-dimensionality assumptions.

Assume that $H$ is a graded connected Hopf algebra. Thus 
$H=\bigoplus_{n\geq 0}H_n$ and $m$ and $\Delta$ are degree-preserving maps.
We are interested in linear endomorphisms of $H$ which preserve
the grading and are zero except on finitely many components:
\[\EndGr(H):=\bigoplus_{n\geq 0}\End(H_n)\,.\]

The following result is central to our constructions.

\begin{proposition}\label{P:interpolation}
The composition, convolution, and smash products of $\End(H)$ restrict to $\EndGr(H)$. Moreover, if $f\in\End(H_p)$, $g\in\End(H_q)$ then
 \begin{equation}\label{E:smash-degree}
f\smashprod g\ \in \bigoplus_{n=\max(p,q)}^{p+q} \End(H_n)
\end{equation}
 and the top and bottom components of $f\smashprod g$ are
\[(f\smashprod g)_{p+q}=f\ast g \text{ \ and, if $p=q$,\ } (f\smashprod g)_p=g\circ f\,.\]
\end{proposition}

Thus the smash product interpolates between  the composition and
convolution products. The analogous interpolation property at all other levels
(permutations, non-commutative symmetric functions, symmetric functions) is a consequence of this general result.

In order to  specialize this construction we let 
\[H=T(V):=\bigoplus_{n\geq 0} V^{\otimes n}\]
 be the tensor algebra of a vector space $V$. It is a graded connected Hopf algebra with coproduct uniquely determined by 
 \[v\mapsto 1\otimes v+v\otimes 1 \text{ \ for \ } v\in V\,.\]
  
 The general linear  group $\GL(V)$ acts  on $V$ and hence on each $V^{\otimes n}$, diagonally. Schur-Weyl duality states that the only endomorphisms of $T(V)$ which commute with the action of $GL(V)$ are (linear combinations of) permutations. Let
 \[\calS:=\bigoplus_{n\geq 0}\field S_n\]
 be the direct sum of all symmetric group algebras.
 
 \begin{lemma}\label{L:SWduality} (Schur-Weyl duality). Suppose $\dim V=\infty$. Then
\[\calS= \EndGr_{GL(V)}\bigl(T(V)\bigr)\,.\]
\end{lemma}
Here each $\sigma\in S_n$ is viewed as an endomorphism of $V^{\otimes n}$ by means of its right action:
\[v_1\cdots v_n\mapsto v_{\sigma(1)}\cdots v_{\sigma(n)}\,.\]

Malvenuto and Reutenauer~ \cite{MR} deduce from here  that $\calS$ is closed under convolution: since the product and coproduct of $T(V)$ commute with the action of
$GL(V)$, the convolution of two permutations must be a linear combination of permutations, by Schur-Weyl duality. The same argument gives us:

\begin{corollary}\label{C:smash-per} The space $\calS$ is closed under the smash product of endomorphisms.
\end{corollary}

This conceptual argument is important because it can  be applied 
to other dualities than Schur-Weyl's, i.e., to centralizer algebras of 
groups  (or even Hopf algebras)  acting on the tensor algebra other than the general linear group. 
It can also be applied to other products of endomorphisms, a remarkable case
being that of {\em Drinfeld} product, which we intend to study in future work.

 Malvenuto and Reutenauer give the following explicit formula for the
 convolution of permutations $\sigma\in\S_p$ and $\tau\in\S_q$:
\begin{equation}\label{E:con-per}
\sigma\ast\tau = \sum_{\xi\in\Sh(p,q)}\xi\circ(\sigma\times\tau)\,,
\end{equation}
where $\Sh(p,q)=\{\xi\in \S_{p+q}\mid \xi(1)<\cdots<\xi(p),\  \xi(p+1)<\cdots<\xi(p+q)\}$ is the set of $(p,q)$-shuffles. Similarly, we find:
\begin{equation}\label{E:smash-per}
\sigma\smashprod\tau = \!\!\!\!\!\sum_{\substack{\max(p,q)\leq n\leq p+q \\
\xi\in\Sh(p,n-p) \\ \eta\in\Sh(p+q-n,n-q)}}
\!\!\!\!\!\xi\circ\bigl((\sigma\circ\eta)\times 1_{n-p}\bigr)\circ
\beta_{2n-p-q,p+q-n}\circ(1_{n-q}\times\tau)
\end{equation}
where $\beta_{u,v}$ is the shuffle of maximum length in $\Sh(u,v)$. 
According to Proposition~\ref{P:interpolation}, we must have
\[(\sigma\smashprod \tau)_{p+q}=\sigma\ast \tau \text{ \ and, if $p=q$,\ } (\sigma\smashprod \tau)_p=\sigma\circ\tau\,,\]
results which may be directly verified from~\eqref{E:con-per} and~\eqref{E:smash-per}. Note that since the action of $S_n$ on $V^{\otimes n}$ is from the right,
composition of permutations corresponds to composition of endomorphisms in the opposite order.

We thus have three associative products on $\calS$: composition (the usual product in each symmetric group algebra $\field S_n$), convolution (which produces an element of degree $p+q$ out of two elements of degrees $p$ and $q$), and the smash product, which produces elements of various degrees
interpolating between the previous two, and is still associative.

\section{The smash product of non-commutative symmetric functions}\label{S:sigma}

The descent set of a permutation $\sigma\in S_n$ is
\[\Des(\sigma):=\{i\in[n{-}1]\,\mid\,\sigma(i)>\sigma(i+1)\}\,.\]
Given $J\subseteq [n{-}1]$, let 
\begin{equation}\label{E:def-X}
X_J:= \sum_{\substack{\sigma\in\S_n \\ \Des(\sigma)\subseteq I}}\!\!\sigma \in \field S_n\,.
\end{equation}
It is convenient to index basis elements of $\Sigma_n$ by compositions of $n$
by means of the bijection 
\[(a_1,a_2,\ldots,a_r)\leftrightarrow \{a_1,a_1+a_2,\ldots,a_1+\cdots+a_{r-1}\}\,.\]
For instance, if $n=9$,
$X_{(1,2,4,2)}=X_{\{1,3,7\}}$.

Let $\Sigma_n$ be the subspace of $\field\S_n$ linearly spanned by the elements $X_\alpha$ as $\alpha$ runs over all compositions of $n$ and
\[\Sigma:=\bigoplus_{n\geq 0}\Sigma_n\,.\]
A fundamental result of Solomon~\cite{Sol} states that $\Sigma_n$ is a subalgebra of the symmetric group algebra $\field\S_n$. $(\Sigma_n,\circ)$ is Solomon's descent algebra. Thus $\Sigma$ is closed under the composition product. It is also well-known that
$\Sigma$ is closed under the convolution product~\cite{GKLLRT,Ges,MR}; in fact,
\begin{equation}\label{E:con-X}
X_{(a_1,\ldots,a_r)}\ast X_{(b_1,\ldots,b_s)}=X_{(a_1,\ldots,a_r,b_1,\ldots,b_s)}\,.
\end{equation}
$(\Sigma,\ast)$ is the algebra of non-commutative symmetric functions.

The following theorem generalizes the two previous results in view of the interpolation property of the smash product.
 
\begin{theorem}\label{T:smash-sigma}
 $\Sigma$ is closed under the smash product.
\end{theorem}

We proceed to sketch four different proofs of this result, each of which is
interesting in its own right.

The most conceptual, but less explicit, proof is based on an important result
of Garsia and Reutenauer~\cite{GReu} which charactertizes
the elements of $\calS$ which belong to $\Sigma$ in terms of the action on the tensor algebra. Let  $L(V)$ be the smallest subspace of $T(V)$ containing $V$ and closed  under $[x,y]:=xy-yx$. $T(V)$ is the free associative algebra on $V$,
$L(V)$ is the free Lie algebra on $V$.
\begin{theorem} ~\cite{GReu}
Let $\varphi\in\calS$. Then $\varphi\in\Sigma$ if and only if
for every  $P_{1},\ldots,P_k \in L(V)$, the subspace spanned by
\[\{P_{\tau(1)}\cdots P_{\tau(k)} \mid \tau\in S_k\}\]
is invariant under the right action of $\varphi$.
\end{theorem}
This implies $\Sigma$ is closed under the smash product. The classical fact that $L(V)$ can also be described as the primitive elements of $T(V)$ allows us to
handle the smash product of endomorphisms~\eqref{E:smash-end} easily, and to
verify that the invariance property is preserved. This proof is valuable because it
can be extended to other situations. For instance, the same argument shows that
$\Sigma$ is closed under the Drinfeld product. 

A more laborious but also more informative proof consists in obtaining an
explicit formula for the smash product of two basis elements of $\Sigma$,
by a direct combinatorial analysis of~\eqref{E:smash-per} when applied to~\eqref{E:def-X}. The result is expressed in terms of the same matrices
as in Theorem~\ref{T:smash-rep}.

\begin{theorem}\label{T:smash-X}
Let $\alpha\comp p$ and
  $\beta\comp q$ be two compositions. Then
\begin{equation}\label{E:smash-X}
X_\alpha\smashprod X_\beta = \sum_{n=\max(p,q)}^{p+q}
\sum_{M\in\M_{\alpha,\beta}^n} X_{c(M)}
\end{equation}
where $c(M)$ is the composition whose parts are the non-zero entries of $M$, read
from left to right and from top to bottom.
\end{theorem}

As an example  we have the following formula
\begin{equation*}
X_{(1^p)} \smashprod X_{(1^q)} = \sum_{n=\operatorname {max}(p,q)}^{p+q} 
\binom{p}{n-q} \binom{q}{n-p} (p+q-n)! \, X_{(1^n)}. 
\end{equation*}  
where $(1^p)$ is the composition of $p$ with $p$ parts equal to $1$. 
The coefficients arise from  the possible ways to fill the
matrix
 \begin{equation*}
    \begin{array}{cccc|c}
      0 & \ast & \cdots & \ast & n-q \\ 
      \ast & \ast & \cdots & \ast & 1 \\
      \vdots & \vdots & \ddots & \vdots & \vdots \\
      \ast & \ast & \cdots & \ast & 1 \\
      \cline{1-4} 
      n-p & 1  & \cdots & \multicolumn{2}{l}{1}  \\
    \end{array}
  \end{equation*}
with row and column sums as prescribed.
%  In the first row, we can only put $1$'s and $0$'s, exactly $n-q$
%  ones and the rest of zeroes.  Similarly for the first column with
%  $n-p$ ones and the rest of zeroes.  All in all we have exactly
%  $\binom{p}{n-q}\times\binom{q}{n-p}$ possibilities as choices for
%  the first row and column. 
%Hence, to fill out the remaining
%spaces of the matrix, 
%we eliminate the rows and columns whose ``headings'' have been filled
%with a $1$. Then, all we have to consider are
%$\bigl(q - (n-p)\bigr) \times \bigl(p - (n -q)\bigr) = (p+q-n)
%\times (p+q-n)$ matrices with exactly one $1$ en each row and
%column, that is, we have to consider all the possible permutation square
%matrices of size $p+q-n$, which are $(p+q-n)!$. 

Similarly, one verifies
\begin{equation*}
X_{(1)}^{\#(n)} = \sum_{k=1}^n S(n,k) X_{(1^k)}\,,
\end{equation*}
where the $S(n,k)$ are the Stirling numbers of the second kind.

%\subsection*{Additional example.} As an another suggestive example of
%the combinatorial rule for the coefficients of the smash product, we
%compute the $n$-th powers of $B_{(1)}$:
%\begin{equation*}
%B_{(1)}^{\#(n)} = \sum_{k=1}^n S(n,k) B_{(1^k)},
%\end{equation*}
%where $S(n,k)$ are the Stirling numbers of second kind. This equation
%can be easily proved by induction on $n$, the inductive step being
%\begin{equation*}
%B_{(1)}^{\#(n)}=B_{(1)}^{\#(n-1)}\smashprod B_{(1)} = \sum_{k=1}^{n-1}
%S(n-1,k)B_{(1^k)}\smashprod B_{(1)}.
%\end{equation*}
%The product $B_{(1^k)}\smashprod B_{(1)}$ is computed by filling the
%following $2\times (k+1)$ matrices
%\setcounter{xx}{0}
%\begin{equation*}
%    \begin{array}{ccccc!{\ver}c}
%      0 & 1 & 0 & \cdots & 1 & k-1 \\ 
%      0 & 0 & 1 & \cdots & 0 & 1 \\ \cline{1-5}
%      0 & 1 & 1 & \cdots & \setcounter{xx}{1} 1 &
%    \end{array},
%    \qquad
%    \setcounter{xx}{0}
%    \begin{array}{cccc!{\ver}c}
%      0 & 1 & \cdots & 1 & k \\ 
%      1 & 0 & \cdots & 0 & 1 \\ \cline{1-4}
%      1 & 1 & \cdots & \setcounter{xx}{1} 1 &
%    \end{array}.
%\end{equation*}
%There are $k$ ways to fill the first matrix and all the possibilities
%yield the composition $(1,\dots,1)$ of $k$. And there is only one way
%to fill the second matrix, which produces the composition
%$(1,\dots,1)$ of $k+1$. Substituting back and reindexing the sums we
%obtain
%\begin{equation*}
%  B_{(1)}^{\#(n)}=\sum_{k=1}^n\bigl[ kS(n-1,k)+S(n-1,k-1)\bigr]
%  B_{(1^k)}=\sum_{k=1}^n S(n,k)B_{(1^k)},
%\end{equation*}
%using a well-known recurrence formula for the Stirling numbers.

\medskip

By the interpolation property of the smash product, Theorem~\ref{T:smash-X} contains as special cases rules for the product in Solomon's descent algebra and for the convolution product of two basis elements of $\Sigma$. One readily verifies that the former is precisely the well-known rule of Garsia, Remmel, Reutenauer, and Solomon as given in~\cite[Proposition 1.1]{GReu}, while the latter is~\eqref{E:con-X}.

This proof is important because it allows us to make the connection with the
smash product of representations of the symmetric group. This point is taken up in Section~\ref{S:sym}

There is a third proof which is analogous to the proof Blessenohl and Laue~\cite{BL} that each $\Sigma_n$ is closed under the composition product. This is of a purely combinatorial nature and is based on a description of descent classes as equivalence classes for a certain relation.

The fourth proof is geometric and is based on an extension of the smash product to the Coxeter complex of the symmetric group (that is, to the faces of the permutahedron). This makes a connection with  recent work of Brown, Mahajan,  Schocker, and others on this aspect of the theory of descent algebras~\cite{Bro,AM,Sch}.

\medskip

There is a coproduct $\Delta$ on the space $\calS$ which is compatible with the convolution product, in the sense that $(\calS,\ast,\Delta)$ is a graded connected Hopf algebra~\cite{AS,MR}. $\Delta$ is {\em not} compatible with the composition product of permutations~\cite[Remarque 5.15]{Mal}. 
The space $\Sigma$ is closed under $\Delta$; in fact
\begin{equation}\label{E:coprod-X}
\Delta(X_{(a_1,\ldots,a_r)})=\sum_{\substack{a_i=b_i+c_i \\ 0\leq b_i,c_i }}
X_{(b_1,\ldots,b_r)^{\widehat{\ }}}\otimes X_{(c_1,\ldots,c_r)^{\widehat{\ }}}\,,
\end{equation}
where ${}^{\widehat{\ }}$ indicates that zero parts are ommited. Moreover, it is known that
$\Delta$ is compatible with {\em both} the convolution product and the composition product 
of $\Sigma$.

\begin{theorem}\label{T:Hopf-sigma}
The space $\Sigma$ endowed with the smash product~\eqref{E:smash-X} and the coproduct~\eqref{E:coprod-X} is a connected Hopf algebra. It is cocommutative but not commutative.
\end{theorem}
To prove compatibility between $\Delta$ and $\smashprod$ on $\Sigma$ we make use of the explicit rule~\eqref{E:smash-X}. The existence of the antipode is automatic in a connected Hopf algebra. An explicit formula is given in Theorem~\ref{T:antipode-Q}.

\begin{theorem}\label{T:iso-sigma}
The map $(\Sigma,\ast,\Delta)\to (\Sigma,\smashprod,\Delta)$ given by
\begin{equation}\label{E:iso-sigma}
X_{(a_1,\ldots,a_r)}\mapsto X_{(a_1)}\smashprod\cdots\smashprod X_{(a_r)}
\end{equation}
is an isomorphism of Hopf algebras (which does not preserve gradings).
\end{theorem}

For this reason $(\Sigma,\smashprod,\Delta)$ may be seen as a deformation of
the Hopf algebra of non-commutative symmetric functions.

\section{The smash product of symmetric functions}\label{S:sym}

In this section we connect the constructions of Sections~\ref{S:rep} and~\ref{S:sigma}.
Consider the map $\phi:\Sigma\to\Lambda$ defined by
\[X_\alpha\mapsto h_\alpha\]
for any composition $\alpha$.
In the original work of Solomon on the descent algebra, it is shown that
the composition product of $\Sigma_n$ corresponds to the internal product of $\Lambda_n$ via $\phi$. It is also known that the convolution product of $\Sigma$ corresponds to the external product of $\Lambda$~\cite{GKLLRT,MR} and that
coproduct~\eqref{E:coprod-X} corresponds to coproduct~\eqref{E:coprod-rep} under $\phi$. 
This generalizes as follows.

\begin{theorem}\label{T:sigma-lambda} The map $\phi:(\Sigma,\smashprod,\Delta)\to(\Lambda,\smashprod,\Delta)$ is a morphism of Hopf algebras.
\end{theorem}

This follows from Theorems~\ref{T:smash-rep} and~\ref{T:smash-X}.
In particular this shows that the smash product of representations~\eqref{E:smash-rep} and the coproduct of representations~\eqref{E:coprod-rep} are compatible.

{}From Theorem~\ref{T:iso-sigma} we deduce that $(\Lambda,\ast,\Delta)$ (the Hopf algebra of symmetric functions) and $(\Lambda,\smashprod,\Delta)$ are isomorphic under a non-degree-preserving isomorphism.

\section{The smash coproduct of quasi-symmetric functions} \label{S:quasi}

Let $\bfX=\{x_1,x_2,\ldots\}$ be a countable set, totally ordered by $x_1<x_2<\cdots$. We say that $\bfX$ is an {\em alphabet}. Let $\field[[\bfX]]$ be the algebra of formal power series on $\bfX$ and $\calQ:=\calQ(\bfX)$ the subspace linearly spanned by the elements
\begin{equation}\label{E:def-M}
M_\alpha:=\sum_{i_1<\dots<i_k}x_{i_1}^{a_1}\cdots x_{i_r}^{a_r}
\end{equation}
as $\alpha=(a_1,\ldots,a_r)$ runs over all compositions of $n$, $n\geq 0$.
$\calQ$ is a graded subalgebra of $\field[[\bfX]]$  known as the algebra of quasi-symmetric functions~\cite{Ges}.  This algebra carries two coproducts
$\Delta_\circ$ and $\Delta_\ast$ which are defined via evaluation of quasi-symmetric functions on alphabets. Let $\bfY$ be another alphabet. We can view the disjoint union $\bfX+\bfY$ and the Cartesian product $\bfX\times\bfY$ as
alphabets as follows: on $\bfX+\bfY$ we keep the ordering among the variables of $\bfX$ and among the variables of $\bfY$, and we require that every variable of $\bfX$ precede every variable of $\bfY$. On $\bfX\times\bfY$ we impose the reverse lexicographic order:
\[(x_h,y_i)\leq (x_j,y_k) \Leftrightarrow y_i<y_k \text{ or } (y_i=y_k \text{ and }x_h< x_j)\,.\]

The coproducts are defined by the formulas
\[\Delta_\circ\bigl(f(\bfX)\bigr):=f(\bfX\times\bfY) \text{ \ and \ }\Delta_\ast\bigl(f(\bfX)\bigr):=f(\bfX+\bfY)\,,\]
together with the identification $\calQ(\bfX,\bfY)\cong \calQ(\bfX)\otimes\calQ(\bfX)$
(separation of variables).

Consider the following pairing between
the homogeneous components of degree $n$ of $\calQ$ and $\Sigma$:
\begin{equation}\label{E:pairing}
\langle M_\alpha, X_\beta\rangle=\delta_{\alpha,\beta}\,.
\end{equation}
It is known~\cite{GKLLRT,Ges,MR} that this pairing identifies the product of quasi-symmetric functions with the coproduct~\eqref{E:coprod-X} of $\Sigma$, 
and the coproducts $\Delta_\circ$ and $\Delta_\ast$ with the composition and
convolution products of $\Sigma$. In other words,
\begin{equation*}
\langle fg, u\rangle =  \langle f\otimes g, \Delta(u) \rangle,\quad
\langle \Delta_{\circ}f, u\otimes v\rangle = \langle f, u\circ v\rangle,\quad
\langle \Delta_{\ast}f, u\otimes v\rangle = \langle f, u\ast v\rangle,
\end{equation*}
for any $f,g\in\calQ$ and $u,v\in\Sigma$. Here we set
\[\langle f\otimes g, u\otimes v\rangle=\langle f,u\rangle\langle g,v\rangle\,.\]
 
Let $\Delta_{\smashprod}$ be the coproduct of $\calQ$ dual to the smash product of $\Sigma$:
\[\langle \Delta_{\smashprod}f, u\otimes v\rangle = \langle f, u\smashprod v\rangle\,.\]
Since the smash product is a sum of terms of various degrees~\eqref{E:smash-degree}, the smash coproduct is a finite sum of the form
\[\Delta(f)=\sum_i f_i\otimes f'_i \text{ \ with \ }
0\leq \deg(f_i) \text{\ and\ }\deg(f'_i)\leq \deg(f)\leq \deg(f_i)+\deg(f'_i)\,.\]
The summands corresponding to $\deg(f)= \deg(f_i)=\deg(f'_i)$ and to
$\deg(f)=\deg(f_i)+\deg(f'_i)$ are the coproducts $\Delta_\circ(f)$ and $\Delta_\ast(f)$, respectively.

\medskip

Let $\bfone+\bfX$ denote the alphabet $\bfX$ together with a new variable $x_0$ smaller than all the others. Let
\[(\bfone+\bfX)\times(\bfone+\bfY)-\bfone\]
 be the Cartesian product of the alphabets $\bfone+\bfX$ and $\bfone+\bfX$ with reverse lexicographic ordering and with the
variable $(x_0,y_0)$ removed.

The following result was obtained in conversation with Arun Ram.

\begin{theorem}\label{T:coprod-Q} For any $f\in\calQ$,
\[\Delta_{\smashprod}\bigl(f(\bfX)\bigr)=f\bigl((\bfone+\bfX)\times(\bfone+\bfY)-\bfone\bigr)\,.\]
\end{theorem}
The set $(\bfone+\bfX)\times(\bfone+\bfY)-\bfone$ can be identified with the disjoint union of $\bfX$, $\bfY$, and $\bfX\times\bfY$, so the evaluation of $f$ on this alphabet produces an element of $\calQ(\bfX,\bfY)\cong\calQ(\bfX)\otimes\calQ(\bfY)$.

\medskip 

Endowed with the coproduct $\Delta_{\smashprod}$, the algebra $\calQ$ is a  Hopf algebra, in duality with the connected Hopf algebra $(\Sigma,\smashprod,\Delta)$ by means of~\eqref{E:pairing}.  We turn to the antipode of this Hopf algebra. 

First, define the evaluation of  quasi-symmetric functions on the
the opposite of an alphabet $\bfX$ by the equation
\begin{equation*}
M_\alpha(-\bfX) = (-1)^k \sum_{i_1\ge \cdots \ge i_k} x^{a_1}_{i_1}
\cdots x^{a_r}_{i_r},
\end{equation*}
for any composition $\alpha=(a_1,\ldots,a_r)$. Compare with~\eqref{E:def-M}.
Second, define the alphabet 
\begin{equation*}
\bfX^* := \bfX+\bfX^2+\bfX^3+\cdots
\end{equation*}
as the disjoint union of the Cartesian powers $\bfX^n$ under reverse lexicographic order. For instance $(x_3,x_1,x_2)<(x_2,x_2)<(x_1,x_3,x_2)$. 

\begin{theorem}\label{T:antipode-Q}
The antipode of the Hopf algebra of quasi-symmetric functions $\calQ$ endowed with  the smash coproduct is
\begin{equation*}
S_{\smashprod}(f) = f(-\bfX^*),
\end{equation*}
\end{theorem}
%This means we evaluate $f$ on the opposite of the alphabet $\bfX^*$.

%\medskip

We define the {\em exponential} of an alphabet $\bfX$ as 
\begin{equation*}
\bfe(\bfX)=\bfX+\bfX^{(2)}+\bfX^{(3)}+\cdots\,,
\end{equation*}
where the {\em divided power} $\bfX^{(n)}$ is the set 
\[\bfX^{(n)}=\{(x_{i_1},x_{i_2},\cdots x_{i_n})\in\bfX^n
\mid x_{i_1}<x_{i_2}<\cdots< x_{i_n} \}\,.\]
 We endow $\bfe(\bfX)$ with the reverse lexicographic order, so that $\bfe(\bfX)$ is a subalphabet of~$\bfX^*$.
Let $\QC:=\prod_{n\geq 0}\calQ_n$.  Given $f\in\QC$, define
\[\varphi(f):=f\bigl(\bfe(\bfX)\bigr)\,.\]
If $f=\sum_{n\geq 0}f_n$ with $f_n$ of degree $n$ then 
 the only the terms in $\varphi(f)$ which contribute to the component of degree $n$ are those $\varphi(f_i)$ for which $i\leq n$. Thus $\varphi(f)$ is a well-defined element of $\QC$.

Under the pairing~\eqref{E:pairing}, $\QC$ identifies with the full linear dual of the space $\Sigma$.
\begin{theorem} The map $\varphi:\QC\to\QC$ is dual to the isomorphism $(\Sigma,\ast,\Delta)\cong(\Sigma,\smashprod,\Delta)$
of Theorem~\ref{T:iso-sigma}.
\end{theorem}

 ÊÊ\bibliography{db}
\bibliographystyle{amsplain}

\end{document}